\documentclass[11pt]{article}
\usepackage[a4paper,margin=26mm]{geometry}
\usepackage{amsmath,amssymb,amsthm,mathtools,microtype,needspace}
\usepackage[T1]{fontenc}
\usepackage{lmodern}
\usepackage[hidelinks]{hyperref}
\hypersetup{pdftitle={Pathwise uniqueness for the stochastic heat equation with a 3/4-H\"older diffusion coefficient},
  pdfauthor={Pei-Sen Li, Jie Xiong, Xu Yang},
  pdfsubject={Research manuscript},pdfkeywords={stochastic heat equation, pathwise uniqueness, critical Holder exponent, white noise, linear comparison}}
\numberwithin{equation}{section}
\newtheorem{theorem}{Theorem}[section]
\newtheorem{definition}[theorem]{Definition}
\newtheorem{condition}{Condition}

\newtheorem{lemma}[theorem]{Lemma}
\newtheorem{proposition}[theorem]{Proposition}

\theoremstyle{remark}

\theoremstyle{plain}
\newcommand{\E}{\mathbb E}
\newcommand{\R}{\mathbb R}
\newcommand{\F}{\mathcal F}
\newcommand{\ind}{\mathbf 1}
\newcommand{\condH}{\textnormal{\textbf{\ref{cond:H}}}}
\newcommand{\e}{\mathrm e}
\newcommand{\fieldnorm}[1]{\lvert\mkern-2mu\lvert\mkern-2mu\lvert#1\rvert\mkern-2mu\rvert\mkern-2mu\rvert_\upsilon}
\renewcommand{\d}{\mathop{}\!\mathrm d}
\title{Pathwise uniqueness for the stochastic heat equation\\
with a $3/4$-H\"older diffusion coefficient}
\author{Pei-Sen Li\thanks{\raggedright School of Mathematics and Statistics,
Beijing Institute of Technology, Beijing, China.
\mbox{E-mail: \href{mailto:peisenli@bit.edu.cn}{\nolinkurl{peisenli@bit.edu.cn}}.}}
\and Jie Xiong\thanks{\raggedright Department of Mathematics and SUSTech International
Center for Mathematics, Southern University of Science and Technology,
Shenzhen, China.
\mbox{E-mail: \href{mailto:xiongj@sustech.edu.cn}{\nolinkurl{xiongj@sustech.edu.cn}}.}}
\and Xu Yang\thanks{\raggedright School of Mathematics and Information Science,
North Minzu University, Yinchuan, China.
\mbox{E-mail: \href{mailto:xuyang@mail.bnu.edu.cn}{\nolinkurl{xuyang@mail.bnu.edu.cn}}.}}}
\date{}
\begin{document}
\makeatletter
\let\standardfnsymbol\@fnsymbol
\renewcommand{\@fnsymbol}[1]{%
  \ifcase#1\or\textdagger\or\textdaggerdbl\or\textasteriskcentered
  \else\standardfnsymbol{#1}\fi}
\maketitle
\let\@fnsymbol\standardfnsymbol
\makeatother

\begin{abstract}
For every $\gamma\in[3/4,1]$, we prove pathwise uniqueness and strong
existence for the stochastic heat equation on $\R$ with space-time white
noise, Lipschitz drift, and a diffusion coefficient locally
$\gamma$-H\"older under suitable growth bounds.
The result settles the critical case $\gamma=3/4$ left open by
Mueller, Mytnik and Perkins.
The proof is based on a newly developed short-time linear comparison
method. On each short interval, we compare the difference of two solutions
with an auxiliary linear process. An expected weighted $L^1$ bound
independent of the noise multiplier, together with a vanishing sum of
approximation errors, yields pathwise uniqueness.
The method may also apply to other SPDEs admitting analogous linear
comparison and stochastic convolution estimates.
\end{abstract}

\section{Introduction and main results}

We consider the stochastic heat equation
\begin{equation}\label{eq:spde}
 \partial_t u_t(x)=\tfrac12\partial_{xx}u_t(x)+b(t,x,u_t(x))
      +\sigma(t,x,u_t(x))\dot W(t,x),\qquad t>0,\ x\in\R,
\end{equation}
where $\dot W$ is space-time Gaussian white noise and
$b,\sigma:\R_+\times\R^2\to\R$ are Borel functions. The initial function
$u_0$ is deterministic and belongs to
\[
 C_{\rm tem}\coloneqq \bigl\{f\in C(\R): |f|_\mu<\infty\text{ for every }\mu>0\bigr\},
 \qquad |f|_\mu\coloneqq \sup_{x\in\R}\e^{-\mu|x|}|f(x)|.
\]
We equip $C_{\rm tem}$ with the topology induced by these seminorms.

For Lipschitz coefficients of linear growth, existence and pathwise
uniqueness follow by Picard iteration; see Walsh
\cite[Chapter~3]{Walsh} on a bounded interval and Shiga
\cite[Theorem~2.2]{Shiga} on $\R$.

Branching particle models give a reason to weaken the Lipschitz condition
on $\sigma$. Konno and Shiga \cite[Theorem~1.4]{KS} derived the equation
for the nonnegative density $u_t$ of a one-dimensional super-Brownian
motion $X_t(\d x)=u_t(x)\,\d x$. With the usual normalization, $u$
satisfies \eqref{eq:spde} with $b=0$ and $\sigma(z)=\sqrt z$ for $z\ge0$.
Exponential duality gives uniqueness in law among nonnegative solutions,
but pathwise uniqueness remains open (see \cite[Example~1]{MMP}).
See also Xiong \cite{Xiong} for a related SPDE formulation and Chen
\cite{Chen} for super-Brownian motion with immigration.

If particles branch at a rate $r(z)$ determined by the local density
$z\ge0$, the noise coefficient becomes $\sqrt{z\,r(z)}$
(see \cite[Section~1]{MP}). For instance, the choice $r(z)=\sqrt z$ gives
$\sigma(z)=z^{3/4}$. Particles can also interact by annihilating on contact. In the model of
Mueller, Mytnik and Perkins \cite[Example~4]{MMP}, each collision
between particles of opposite sign removes equal positive and negative
mass. If $u^+$ and $u^-$ are the two densities in a scaling limit,
their difference $u=u^+-u^-$ can take either sign and satisfies
\eqref{eq:spde} with $b=0$ and $\sigma(z)=\sqrt{|z|}$ for $z\in\R$.
These models lead us to study H\"older continuous diffusion coefficients.
It is natural to ask the following question.

\smallskip\noindent
\emph{How far can the regularity of $\sigma$ be weakened while retaining
pathwise uniqueness among real-valued solutions?}
\par\smallskip

For H\"older diffusion coefficients, this problem remained
open only at the critical exponent $\gamma=3/4$.
Mytnik and Perkins \cite[Theorem~1.2]{MP} proved pathwise uniqueness for
\eqref{eq:spde} when $\sigma$ is locally H\"older continuous of order
$\gamma>3/4$, under their growth hypotheses.
For every $0<\gamma<3/4$, Mueller, Mytnik and Perkins
\cite[Theorem~1.1]{MMP} proved nonuniqueness in law and pathwise
nonuniqueness by constructing a nontrivial real-valued solution with
$b=0$, $\sigma(z)=|z|^\gamma$ and zero initial data.

Our main aim is to settle the critical case $\gamma=3/4$.
We give a concise proof for all $\gamma\in[3/4,1]$ using the short-time
linear comparison method introduced in Subsection~\ref{sec:strategy}.

\begin{samepage}
\begin{condition}\label{cond:H}
For some $\gamma\in[3/4,1]$ and every $T<\infty$, there exist constants
$C,L,K,a,m\ge0$ such that the following hold for $0\le t\le T$
and $x,z,w\in\R$.
\begin{enumerate}
\item[(i)] The coefficients have linear growth:
\[
 |b(t,x,z)|+|\sigma(t,x,z)|\le C(1+|z|).
\]
\item[(ii)] The drift is Lipschitz:
\[
 |b(t,x,z)-b(t,x,w)|\le L|z-w|.
\]
\item[(iii)] For $|z-w|\le1$, the diffusion coefficient satisfies
\[
 |\sigma(t,x,z)-\sigma(t,x,w)|
 \le K\e^{a|x|}(1+|z|+|w|)^m |z-w|^\gamma.
\]
\end{enumerate}
\end{condition}
\end{samepage}

\Needspace{13\baselineskip}
All filtered probability spaces
$(\Omega,\F,(\F_t)_{t\ge0},\mathbb P)$ considered below satisfy the usual conditions.
For a space-time Gaussian white noise $W$ with intensity $\d t\,\d x$
on such a space, we require that the noise up to time $t$ is $\F_t$-measurable and its future
increments are independent of $\F_t$. The natural filtration of $W$ is understood
to be completed and right-continuous.

On a measure space $(E,\nu)$ we identify measurable functions $f:E\to\R$
that agree $\nu$-almost everywhere. For $1\le p<\infty$, set
\[
 L^p(E,\nu)\coloneqq \bigl\{f:\|f\|_{L^p(E,\nu)}<\infty\bigr\},\qquad
 \|f\|_{L^p(E,\nu)}\coloneqq \left(\int_E |f|^p\d\nu\right)^{1/p}.
\]
For $p=\infty$ the norm is the $\nu$-essential supremum of $|f|$.
For $\lambda>0$, set $\rho(x)\coloneqq\e^{-\lambda\sqrt{1+x^2}}$.
On $\R$ we write $\|f\|_p\coloneqq \|f\|_{L^p(\R,\d x)}$ and
$\|f\|_{p,\rho}\coloneqq \|f\|_{L^p(\R,\rho(x)\d x)}$.
For $fg\in L^1(\R,\d x)$ we write $\langle f,g\rangle\coloneqq \int_\R f(x)g(x)\d x$.

\begin{definition}[Weak and strong solutions]\label{def:solution}
Let $u_0\in C_{\rm tem}$ be deterministic, and let
$u=(u_t(x))_{t\ge0,\,x\in\R}$ be a real-valued predictable random field
on some filtered probability space.
\begin{enumerate}
\item[(i)] We call $u$ a \emph{weak solution} of \eqref{eq:spde} with
initial function $u_0$ if $(t\mapsto u_t)\in C(\R_+;C_{\rm tem})$
almost surely and there exists a space-time Gaussian white noise $W$
on the same probability space such that, for each $\phi\in C_c^\infty(\R)$,
all integrals below are well defined and the following identity holds
almost surely for all $t\ge0$:
\[
 \begin{aligned}
 \langle u_t,\phi\rangle={}&\langle u_0,\phi\rangle+
 \int_0^t\!\bigl[\tfrac12\langle u_s,\phi''\rangle+
       \langle b(s,\cdot,u_s),\phi\rangle\bigr]\d s\\
 &+\int_0^t\!\int_\R\phi(x)\sigma(s,x,u_s(x))W(\d s,\d x).
 \end{aligned}
\]
\item[(ii)] A weak solution $u$ driven by $W$ is called a \emph{strong solution}
if $u$ is adapted to the natural filtration of $W$.
\end{enumerate}
Unless otherwise specified, a \emph{solution} of \eqref{eq:spde} means a weak solution.
We say that \emph{pathwise uniqueness} holds if any two solutions on the same
filtered probability space, driven by the same white noise and with the
same initial function, are indistinguishable.
\end{definition}

\begin{theorem}\label{thm:main}
Suppose that \condH{} holds. Then, for every deterministic $u_0\in C_{\rm tem}$,
equation \eqref{eq:spde} has a pathwise unique strong solution with initial
function $u_0$.
\end{theorem}

The joint law of $(u,W)$ is therefore uniquely determined by $u_0$. In particular, consider
\[
 \partial_t u=\tfrac12\partial_{xx}u+|u|^{3/4}\dot W
\]
with zero initial function. Since
$\big||z|^{3/4}-|w|^{3/4}\big|\le\big||z|-|w|\big|^{3/4}\le|z-w|^{3/4}$
and $u\equiv0$ is a solution, Theorem~\ref{thm:main} gives $u\equiv0$
almost surely for every real-valued solution.

It is worth noting that we impose no nondegeneracy condition on $\sigma$.
This is important for branching models, where the noise vanishes at zero
population density. Under uniform nondegeneracy, Han \cite{HanErgodicity}
proved uniqueness in law for stochastic heat equations on bounded
intervals with locally H\"older drift and $\gamma$-H\"older diffusion,
$\gamma>3/4$. His recent work
\cite[Theorem~1.1]{Han} extends uniqueness in law to every positive global
H\"older order for bounded, uniformly nondegenerate diffusion in the
drift-free scalar equation on the one-dimensional torus.

\subsection{The short-time linear comparison method}\label{sec:strategy}
In this subsection, we briefly describe our method for proving pathwise uniqueness.
Let $u,v$ be two solutions driven by the same white noise $W$.
At each time $r$ and position $x$, define the solution difference $D_r(x)$,
the noise coefficient difference $A_r(x)$, and the drift difference $G_r(x)$ by
\[
 \begin{aligned}
 D_r(x)&\coloneqq u_r(x)-v_r(x),\\
 A_r(x)&\coloneqq \sigma(r,x,u_r(x))-\sigma(r,x,v_r(x)),\\
 G_r(x)&\coloneqq b(r,x,u_r(x))-b(r,x,v_r(x)).
 \end{aligned}
\]
Fix $\lambda>0$. For equal initial functions, our goal is to show that
$\E\|D_t\|_{1,\rho}=0$ for every $t$.

Mytnik and Perkins \cite{MP} use the Yamada--Watanabe method \cite{YW}
to prove pathwise uniqueness for $\gamma>3/4$.
To control the local-time term arising from It\^o's formula, they obtain
bounds of order $a_n^{\gamma-3/4}$, where $a_n\downarrow0$
(see \cite[(3.18)]{MP}). At $\gamma=3/4$, this factor provides no decay.

To treat the endpoint $\gamma=3/4$, we introduce the \emph{short-time
linear comparison method}.
Fix $t>0$, put $h=t/N$ and $t_j=jh$.
On each interval $[t_j,t_{j+1}]$, we introduce the following auxiliary
linear equation for the \emph{linear comparison process} $Z^j$.
\[
 \partial_r Z^j=\tfrac12\partial_{xx}Z^j+G+\beta^{(h)} Z^j\dot W,
 \qquad Z^j_{t_j}=D_{t_j}.
\]
Here $\beta^{(h)}$ is the quotient $A/D$ (zero where $D=0$) truncated to
$[-\tfrac12h^{-1/4},\tfrac12h^{-1/4}]$.
We compare $Z^j$ with the difference $D=u-v$, which satisfies
\[
 \partial_r D=\tfrac12\partial_{xx}D+G+A\dot W .
\]
Subtracting the equation for $Z^j$ cancels the drift $G$ and gives
\[
 \partial_r(D-Z^j)=\tfrac12\partial_{xx}(D-Z^j)
       +\bigl[\beta^{(h)}(D-Z^j)+(A-\beta^{(h)}D)\bigr]\dot W ,
 \qquad (D-Z^j)_{t_j}=0 .
\]
We have
\[
 \E\|D_{t_{j+1}}\|_{1,\rho}
 \le \E\|Z^j_{t_{j+1}}\|_{1,\rho}
    +\underbrace{\E\|D_{t_{j+1}}-Z^j_{t_{j+1}}\|_{1,\rho}}
                _{\varepsilon_j}.
\]

\smallskip\noindent\emph{Step 1.}
We first estimate $\E\|Z^j_{t_{j+1}}\|_{1,\rho}$ in terms of the difference
at $t_j$ and the drift on $[t_j,t_{j+1}]$.
As $h\downarrow0$, the stochastic integral is taken over a shorter
time interval, while $|\beta^{(h)}|$ may grow.
Therefore we need a precise estimate that controls both effects.
Proposition~\ref{prop:linear-comparison}(ii) gives
\begin{equation}\label{eq:method-linear-estimate}
 \E\|Z^j_{t_{j+1}}\|_{1,\rho}
 \le \e^{ch}\E\|D_{t_j}\|_{1,\rho}
     +L\int_{t_j}^{t_{j+1}}\e^{c(t_{j+1}-s)}\E\|D_s\|_{1,\rho}\,\d s .
\end{equation}
Here $c=\lambda^2/2$ depends only on the weight $\rho$, and $L$ is the Lipschitz
constant of the drift. Both are independent of $h$ and $\beta^{(h)}$.

\smallskip\noindent\emph{Step 2.}
We next show that the sum of the errors $\varepsilon_j$ tends to zero as $h\downarrow0$.
Under the bound $|A|\le K|D|^\gamma$ with $\gamma\ge3/4$, summing the errors over the $t/h$
intervals as in \eqref{eq:sum} gives
\begin{equation}\label{eq:method-error-sum}
 \sum_{j=0}^{N-1}\varepsilon_j
 =\left(\frac{t}{h}\right)^{3/4}o(h^\gamma)
 =o\bigl(h^{\gamma-3/4}\bigr)\longrightarrow0
 \qquad\text{as }h\downarrow0 .
\end{equation}
\smallskip\noindent
Substituting \eqref{eq:method-linear-estimate} into the bound for
$\E\|D_{t_{j+1}}\|_{1,\rho}$, multiplying by $\e^{-ct_{j+1}}$, and summing
over $j$, we obtain, with $Y(s)=\e^{-cs}\E\|D_s\|_{1,\rho}$,
\[
 Y(t)\le Y(0)+L\int_0^t Y(s)\,\d s+\sum_{j=0}^{N-1}\e^{-ct_{j+1}}\varepsilon_j.
\]
The sum tends to zero by \eqref{eq:method-error-sum}. Since $Y(0)=0$,
Gronwall's inequality gives $Y(t)=0$ for every $t$.

\section{The linear comparison estimate}\label{sec:linear}

Let $W$ be a space-time Gaussian white noise on the given filtered probability
space with filtration $(\F_r)_{r\ge0}$. Let $(P_r)_{r\ge0}$ denote the heat
semigroup on $\R$ for $\tfrac12\partial_{xx}$. Thus $P_0 f=f$ and, for $r>0$,
\[
 (P_r f)(x)=\int_\R p_r(x-y)f(y)\d y,\qquad
 p_r(z)\coloneqq (2\pi r)^{-1/2}\e^{-z^2/(2r)} .
\]
Fix $\lambda>0$ and put $c\coloneqq\lambda^2/2$.
Throughout this section, equalities and inequalities
between random functions are understood in $L^2(\R,\rho(x)\d x)$ almost surely
at each fixed time.

Fix $0\le s<t$. Let $\xi$ be an $\F_s$-measurable
$L^2(\R,\rho(x)\d x)$-valued random variable, and let $G,\beta$ be
predictable real-valued random fields on $(s,t]\times\R$.
Consider the linear equation
\begin{equation}\label{eq:linear-white}
 \begin{aligned}
 Z_r(x)
 &=(P_{r-s}\xi)(x)+\int_s^r(P_{r-q}G_q)(x)\,\d q\\
 &\quad+\int_s^r\!\int_\R p_{r-q}(x-y)\beta_q(y)Z_q(y)
                    \,W(\d q,\d y),\qquad s\le r\le t.
 \end{aligned}
\end{equation}

\begin{proposition}[Linear comparison estimate]\label{prop:linear-comparison}
Assume that $\E\|\xi\|_{2,\rho}^2<\infty$,
$\int_s^t\E\|G_q\|_{2,\rho}^2\d q<\infty$ and $|\beta|\le K_\beta$
for some deterministic constant $K_\beta<\infty$.
Then \eqref{eq:linear-white} has a unique predictable solution $Z$ in the class
$\sup_{s\le r\le t}\E\|Z_r\|_{2,\rho}^2<\infty$.
\begin{enumerate}
\item[(i)] If $\xi$ and $G$ are nonnegative, then $Z_r\ge0$ for each $r\in[s,t]$.
\item[(ii)] For signed data and each $r\in[s,t]$, we have
\begin{equation}\label{eq:source-mass}
 \E\|Z_r\|_{1,\rho}
 \le \e^{c(r-s)}\E\|\xi\|_{1,\rho}
                 +\int_s^r\e^{c(r-q)}\E\|G_q\|_{1,\rho}\d q .
\end{equation}
\end{enumerate}
\end{proposition}

To prove \eqref{eq:source-mass}, we first establish nonnegativity for the equation
with spatially smoothed noise and pass to the white-noise limit.
For signed data, linearity then bounds $|Z|$ by the solution with initial
value $|\xi|$ and drift $|G|$. Taking expectations in its mild equation
removes the stochastic integral and gives \eqref{eq:source-mass}.

\begin{lemma}\label{lem:heat-bounds}
For $q\ge0$ and $f\in L^2(\R,\rho(x)\d x)$, we have
\[
 P_q\rho\le\e^{cq}\rho,\qquad
 \|P_q f\|_{2,\rho}^2\le\e^{cq}\|f\|_{2,\rho}^2.
\]
\end{lemma}
\begin{proof}
Differentiating $p_q(z)$ for $q>0$, we obtain
\begin{equation}\label{eq:heat-kernel-pde}
 \partial_qp_q(z)
 =\left(\frac{z^2}{2q^2}-\frac1{2q}\right)p_q(z)
 =\tfrac12\partial_{zz}p_q(z).
\end{equation}
Using integration by parts and \eqref{eq:heat-kernel-pde}, we obtain
\begin{equation}\label{eq:heat-weight-derivative}
 \begin{aligned}
 \partial_q(P_q\rho)(x)
 &=\int_\R\partial_qp_q(x-y)\rho(y)\d y\\
 &=\tfrac12\int_\R\partial_{yy}p_q(x-y)\rho(y)\d y
 =\tfrac12(P_q\rho'')(x).
 \end{aligned}
\end{equation}
Differentiating $\rho$ twice, we obtain
\begin{equation}\label{eq:rho-generator-bound}
 \frac{\rho''(x)}{2\rho(x)}
 =\frac{\lambda^2x^2}{2(1+x^2)}
  -\frac{\lambda}{2(1+x^2)^{3/2}}\le c.
\end{equation}
By \eqref{eq:heat-weight-derivative}, \eqref{eq:rho-generator-bound}
and positivity of $P_q$, we obtain
\[
 \partial_q\bigl(\e^{-cq}P_q\rho\bigr)
 =\e^{-cq}P_q\bigl(\tfrac12\rho''-c\rho\bigr)\le0.
\]
By monotonicity and $\lim_{q\downarrow0}P_q\rho=\rho$, we obtain the first bound.
Finally, applying Jensen's inequality and the first bound, we obtain
\[
 \|P_q f\|_{2,\rho}^2
 \le\int_\R\rho(x)P_q(|f|^2)(x)\d x
 =\langle P_q\rho,|f|^2\rangle
 \le\e^{cq}\|f\|_{2,\rho}^2,
\]
which is the second bound.
\end{proof}

\begin{lemma}\label{lem:stochastic-convolution}
Fix $0\le s<r$ and let $B=(B_q(y))_{s<q<r,\,y\in\R}$ be a real-valued
predictable random field satisfying
\begin{equation}\label{eq:convolution-integrability}
 \int_s^r \frac{\E\|B_q\|_{2,\rho}^2}{2\sqrt{\pi(r-q)}}\,\d q<\infty.
\end{equation}
Set $H_{q,x}(y)\coloneqq p_{r-q}(x-y)B_q(y)$. Then the integral below is
well defined in $L^2(\R,\rho(x)\d x)$, and we have
\begin{equation}\label{eq:isometry}
 \begin{aligned}
 \E\left\|\int_s^r\!\int_\R H_{q,\cdot}(y)\,W(\d q,\d y)\right\|_{2,\rho}^2
 &=\int_\R\rho(x)\,\E\int_s^r\|H_{q,x}\|_2^2\,\d q\,\d x\\
 &\le\int_s^r \frac{\e^{c(r-q)/2}}{2\sqrt{\pi(r-q)}}\E\|B_q\|_{2,\rho}^2\,\d q .
 \end{aligned}
\end{equation}
\end{lemma}
\begin{proof}
Using $p_\ell^2=p_{\ell/2}/(2\sqrt{\pi \ell})$ and Lemma~\ref{lem:heat-bounds}, we obtain
\begin{equation}\label{eq:weighted-kernel-square}
 \int_\R\rho(x)p_\ell(x-y)^2\d x
 =\frac{(P_{\ell/2}\rho)(y)}{2\sqrt{\pi \ell}}
 \le\frac{\e^{c\ell/2}\rho(y)}{2\sqrt{\pi \ell}},\qquad \ell>0.
\end{equation}
Taking $\ell=r-q$ in \eqref{eq:weighted-kernel-square} gives
\[
 \begin{aligned}
 &\int_s^r\!\int_\R \E|B_q(y)|^2
   \left(\int_\R\rho(x)p_{r-q}(x-y)^2\d x\right)\d y\,\d q\\
 &\qquad\le\int_s^r \frac{\e^{c(r-q)/2}}{2\sqrt{\pi(r-q)}}\E\|B_q\|_{2,\rho}^2\,\d q<\infty,
 \end{aligned}
\]
where the last integral is finite by \eqref{eq:convolution-integrability}
and $\e^{c(r-q)/2}\le\e^{c(r-s)/2}$.
Since $\rho>0$, we have $\E\int_s^r\|H_{q,x}\|_2^2\,\d q<\infty$ for
almost every $x$. The white-noise isometry now gives the equality
in \eqref{eq:isometry}.
\end{proof}

\begin{lemma}\label{lem:smoothed-noise-field}
Fix $\varepsilon>0$. Then the random field
\[
 (r,z)\longmapsto\int_0^r\!\int_\R p_\varepsilon(z-y)\,W(\d q,\d y),
 \qquad (r,z)\in\R_+\times\R,
\]
admits a jointly continuous version, which we denote by
$(W^\varepsilon_r(z))_{r\ge0,\,z\in\R}$. For each fixed $z$,
$(W^\varepsilon_r(z))_{r\ge0}$ is an $(\F_r)_{r\ge0}$-Brownian motion
with variance rate $p_{2\varepsilon}(0)$.
\end{lemma}

\begin{proof}
Every finite linear combination of the values of the field
is a white-noise integral with a deterministic integrand and hence is centered
Gaussian. By the white-noise isometry and the heat semigroup property, we obtain
\[
 \E[W^\varepsilon_r(x)W^\varepsilon_q(z)]
 =(r\wedge q)\int_\R p_\varepsilon(x-y)p_\varepsilon(z-y)\,\d y
 =(r\wedge q)p_{2\varepsilon}(x-z).
\]
For $r,q\le T$, this gives
\[
 \begin{aligned}
 \E|W^\varepsilon_r(x)-W^\varepsilon_q(z)|^2
 &=|r-q|p_{2\varepsilon}(0)
   +2(r\wedge q)\bigl(p_{2\varepsilon}(0)-p_{2\varepsilon}(x-z)\bigr)\\
 &\le C_{\varepsilon,T}\bigl(|r-q|+|x-z|^2\bigr).
 \end{aligned}
\]
Using Gaussian sixth moments and Kolmogorov's continuity theorem in $(r,x)$,
we obtain the stated jointly continuous version. For fixed $z$ and $q<r$,
the increment $W^\varepsilon_r(z)-W^\varepsilon_q(z)$ uses only the white noise
on $(q,r]\times\R$. It is independent of $\F_q$ and has law
$N(0,(r-q)p_{2\varepsilon}(0))$. Since $W^\varepsilon_0(z)=0$, this proves
the Brownian motion assertion.
\end{proof}

We use martingale measures in the sense of Walsh
\cite[Chapter~2, pp.~287--291]{Walsh}. For each bounded Borel set $A\subset\R$,
$M(A)=(M_r(A))_{r\ge0}$ is a square-integrable martingale starting at zero.
At each time $r$, the map $A\mapsto M_r(A)$ is countably additive in
$L^2(\Omega)$. We call $M$ continuous if each $M(A)$ has continuous paths.
For a rectangle $D=A\times B\times(a,b]$, write
$Q(D)\coloneqq \langle M(A),M(B)\rangle_b-\langle M(A),M(B)\rangle_a$.
The martingale measure is \emph{worthy} if $|Q(D)|\le K(D)$ on rectangles
for a random measure $K$ that is symmetric and positive definite in the
space variables and has finite expected mass on bounded sets, with
$r\mapsto K(A\times B\times(0,r])$ predictable for each $A,B$.

For $H_q(z)=\eta\ind_{(a,b]}(q)\ind_A(z)$ with $\eta$ bounded and
$\F_a$-measurable, the Walsh integral is $\eta(M_b(A)-M_a(A))$.
For a predictable field $H$ on $(0,T]\times\R$, set
\[
 \|H\|_M^2\coloneqq \E\int_{\R^2\times(0,T]} |H_q(x)H_q(y)|\,K(\d x,\d y,\d q).
\]
When $\|H\|_M<\infty$, we approximate $H$ in this norm by finite sums of
elementary fields and take the limit of their integrals in $L^2(\Omega)$
(see \cite[Chapter~2, pp.~292--295]{Walsh}).

For $\varepsilon\ge0$, $r\ge0$ and bounded Borel sets $A\subset\R$, define
\[
 M^\varepsilon_r(A)\coloneqq \int_0^r\!\int_\R(P_\varepsilon\ind_A)(y)\,W(\d q,\d y).
\]
For $\varepsilon>0$ this gives the spatial convolution of $W$ with $p_\varepsilon$.
Since $P_0$ is the identity, $M^0$ is the original white noise $W$.
In particular, for $0\le a<b$ we have
\[
 M^0_b(A)-M^0_a(A)=W((a,b]\times A).
\]

\begin{lemma}\label{lem:smoothed-noise-integral}
For every $\varepsilon\ge0$, $M^\varepsilon$ is a continuous worthy martingale
measure. If $0\le s<r<\infty$ and $H$ is predictable with
$\E\int_s^r\|H_q\|_2^2\,\d q<\infty$, then both integrals below are
square integrable and
\begin{equation}\label{eq:smoothed-noise-integral}
 \int_s^r\!\int_\R H_q(z)\,M^\varepsilon(\d q,\d z)
 =\int_s^r\!\int_\R (P_\varepsilon H_q)(y)\,W(\d q,\d y)
 \qquad\text{a.s.}
\end{equation}
\end{lemma}

\begin{proof}
For $\varepsilon=0$, the white-noise martingales $M^0(A)$ have continuous
paths and covariations $r|A\cap B|$. Their
covariance measure is supported on the spatial diagonal and serves as a
dominating measure. Identity \eqref{eq:smoothed-noise-integral} follows
from $P_0H_q=H_q$. Now fix $\varepsilon>0$.
For a bounded Borel set $A$, the stochastic Fubini theorem
\cite[Theorem~2.6]{Walsh} gives
\begin{equation}\label{eq:smoothed-noise-set}
 \begin{aligned}
 \int_A W^\varepsilon_r(z)\,\d z
 &=\int_A\left(\int_0^r\!\int_\R
        p_\varepsilon(z-y)\,W(\d q,\d y)\right)\d z\\
 &=\int_0^r\!\int_\R
       \left(\int_A p_\varepsilon(z-y)\,\d z\right)W(\d q,\d y)\\
 &=\int_0^r\!\int_\R (P_\varepsilon\ind_A)(y)\,W(\d q,\d y)=M^\varepsilon_r(A).
 \end{aligned}
\end{equation}
The integrability condition in \cite[Theorem~2.6]{Walsh} is
\[
 \int_A\!\int_0^r\!\int_\R p_\varepsilon(z-y)^2\,\d y\,\d q\,\d z
 =r|A|p_{2\varepsilon}(0)<\infty.
\]
By its white-noise representation, each $M^\varepsilon(A)$ is a continuous
square-integrable martingale. For bounded Borel sets $A,B$,
\[
 \langle M^\varepsilon(A),M^\varepsilon(B)\rangle_r
 =r\langle P_\varepsilon\ind_A,P_\varepsilon\ind_B\rangle
 =r\int_A\!\int_B p_{2\varepsilon}(z-y)\,\d y\,\d z .
\]
The bound $\E|M^\varepsilon_r(A)|^2=r\|P_\varepsilon\ind_A\|_2^2\le r|A|$
gives countable additivity in $L^2(\Omega)$ on bounded sets.
The covariance measure $p_{2\varepsilon}(z-y)\,\d z\,\d y\,\d q$
also serves as the dominating measure, so $M^\varepsilon$ is worthy.

For $H_q(z)=\eta\ind_{(a,b]}(q)\ind_A(z)$, with $s\le a<b\le r$ and
$\eta$ bounded and $\F_a$-measurable, both sides of
\eqref{eq:smoothed-noise-integral} equal
\[
 \eta\bigl(M^\varepsilon_b(A)-M^\varepsilon_a(A)\bigr)
 =\eta\int_a^b\!\int_\R (P_\varepsilon\ind_A)(y)\,W(\d q,\d y).
\]
For general $H$, extended by zero outside $(s,r]$, the covariance kernel gives
\[
 \|H\|_{M^\varepsilon}^2
 =\E\int_s^r\|P_\varepsilon |H_q|\|_2^2\,\d q
 \le \E\int_s^r\|H_q\|_2^2\,\d q .
\]
The white-noise isometry and $L^2$-contractivity of $P_\varepsilon$ give
the same upper bound for the second moment of the right-hand side of
\eqref{eq:smoothed-noise-integral}. Approximating $H$ in predictable $L^2$
by finite sums of the elementary fields above proves the identity.
\end{proof}

We now fix $s,t,\xi,G$ and $\beta$ as in the hypotheses of
Proposition~\ref{prop:linear-comparison}. For $\varepsilon\ge0$, we consider
the mild equation
\begin{equation}\label{eq:linear}
 \begin{aligned}
 Z^\varepsilon_r(x)
 &=(P_{r-s}\xi)(x)+\int_s^r(P_{r-q}G_q)(x)\,\d q\\
 &\quad+\int_s^r\!\int_\R p_{r-q}(x-y)\beta_q(y)Z^\varepsilon_q(y)
                    \,M^\varepsilon(\d q,\d y),\qquad s\le r\le t.
 \end{aligned}
\end{equation}
The stochastic integral in \eqref{eq:linear} is the Walsh integral of
Lemma~\ref{lem:smoothed-noise-integral}.

\begin{lemma}\label{lem:linear-approximation}
\begin{enumerate}
\item[(i)] For each $\varepsilon\ge0$, \eqref{eq:linear} has a unique predictable
solution $Z^\varepsilon$ in the class
$\sup_{s\le r\le t}\E\|Z^\varepsilon_r\|_{2,\rho}^2<\infty$. Moreover,
\[
 \sup_{\varepsilon\ge0}\sup_{s\le r\le t}\E\|Z^\varepsilon_r\|_{2,\rho}^2<\infty .
\]
\item[(ii)] For each fixed $r\in[s,t]$, we have
\begin{equation}\label{eq:white}
 \lim_{\varepsilon\downarrow0}\E\|Z^\varepsilon_r-Z^0_r\|_{2,\rho}^2=0 .
\end{equation}
\end{enumerate}
\end{lemma}

\begin{proof}
\emph{(i)} We construct the solution by a contraction argument.
For $\upsilon>0$ and predictable real-valued random fields $Y$ on
$[s,t]\times\R$, set
\[
 \fieldnorm{Y}\coloneqq \sup_{s\le r\le t}\e^{-\upsilon(r-s)/2}
        \big(\E\|Y_r\|_{2,\rho}^2\big)^{1/2},\qquad
 \mathcal X_{s,t}\coloneqq \{Y:\fieldnorm{Y}<\infty\}.
\]
Since $[s,t]$ is bounded, $\fieldnorm{Y}<\infty$ is equivalent to
$\sup_{s\le r\le t}\E\|Y_r\|_{2,\rho}^2<\infty$, so
$\mathcal X_{s,t}$ does not depend on $\upsilon$.
Write \eqref{eq:linear} as $Z^\varepsilon=I+\mathcal T_\varepsilon Z^\varepsilon$, where
\[
 \begin{aligned}
 I_r&\coloneqq P_{r-s}\xi+\int_s^rP_{r-q}G_q\,\d q,\qquad\text{and}\\
 (\mathcal T_\varepsilon Y)_r(x)
 &\coloneqq \int_s^r\!\int_\R p_{r-q}(x-y)\beta_q(y)Y_q(y)\,M^\varepsilon(\d q,\d y).
 \end{aligned}
\]
By Jensen's inequality and Lemma~\ref{lem:heat-bounds}, we have
\[
 \begin{aligned}
 \sup_{s\le r\le t}\E\|I_r\|_{2,\rho}^2
 &\le 2\sup_{s\le r\le t}\left(
     \E\|P_{r-s}\xi\|_{2,\rho}^2
     +(r-s)\int_s^r\E\|P_{r-q}G_q\|_{2,\rho}^2\d q\right)\\
 &\le 2\sup_{s\le r\le t}\left(
     \e^{c(r-s)}\E\|\xi\|_{2,\rho}^2
     +(r-s)\int_s^r\e^{c(r-q)}\E\|G_q\|_{2,\rho}^2\d q\right)\\
 &\le 2\e^{c(t-s)}
     \left(\E\|\xi\|_{2,\rho}^2
     +(t-s)\int_s^t\E\|G_q\|_{2,\rho}^2\d q\right).
 \end{aligned}
\]
Hence $I\in\mathcal X_{s,t}$.
For $Y\in\mathcal X_{s,t}$ and $s<r\le t$, set
$H_{q,x}(y)\coloneqq p_{r-q}(x-y)\beta_q(y)Y_q(y)$. Applying \eqref{eq:isometry}
with $B=\beta Y$ and using $\int_s^r \frac{\d q}{2\sqrt{\pi(r-q)}}=\sqrt{(r-s)/\pi}$, we obtain
\begin{equation}\label{eq:smoothed-integrand-bound}
 \int_\R\rho(x)\,\E\int_s^r\|H_{q,x}\|_2^2\,\d q\,\d x
 \le K_\beta^2\e^{c(r-s)/2}\sqrt{\frac{r-s}{\pi}}\,
       \sup_{s\le q\le r}\E\|Y_q\|_{2,\rho}^2<\infty .
\end{equation}
Since $\rho>0$, for almost every $x$ we may apply
\eqref{eq:smoothed-noise-integral} to $H_{\cdot,x}$. By the white-noise
isometry, $L^2$-contractivity of $P_\varepsilon$ and \eqref{eq:isometry}, we obtain
\begin{equation}\label{eq:linear-iteration}
 \begin{aligned}
 \E\|(\mathcal T_\varepsilon Y)_r\|_{2,\rho}^2
 &=\int_\R\rho(x)\,\E\int_s^r\|P_\varepsilon H_{q,x}\|_2^2\,\d q\,\d x\\
 &\le\int_\R\rho(x)\,\E\int_s^r\|H_{q,x}\|_2^2\,\d q\,\d x\\
 &\le K_\beta^2\int_s^r \frac{\e^{c(r-q)/2}}{2\sqrt{\pi(r-q)}}\E\|Y_q\|_{2,\rho}^2\,\d q .
 \end{aligned}
\end{equation}
Choose $\upsilon>c/2+4K_\beta^4$.
Multiplying \eqref{eq:linear-iteration} by $\e^{-\upsilon(r-s)}$ and using
$\E\|Y_q\|_{2,\rho}^2\le\e^{\upsilon(q-s)}\fieldnorm{Y}^2$, we obtain
\[
 \begin{aligned}
 \fieldnorm{\mathcal T_\varepsilon Y}^2
 &\le K_\beta^2\sup_{s\le r\le t}\e^{-\upsilon(r-s)}
       \int_s^r \frac{\e^{c(r-q)/2}}{2\sqrt{\pi(r-q)}}\E\|Y_q\|_{2,\rho}^2\,\d q\\
 &\le K_\beta^2\fieldnorm{Y}^2\sup_{s\le r\le t}\int_0^{r-s}
       \frac{\e^{-(\upsilon-c/2)\ell}}{2\sqrt{\pi\ell}}\,\d\ell\\
 &\le K_\beta^2 \fieldnorm{Y}^2\int_0^\infty\frac{\e^{-(\upsilon-c/2)\ell}}{2\sqrt{\pi\ell}}\,\d\ell\\
 &=\frac{K_\beta^2}{2\sqrt{\upsilon-c/2}}\fieldnorm{Y}^2
 \le\tfrac14\fieldnorm{Y}^2.
 \end{aligned}
\]
For $Y_1,Y_2\in\mathcal X_{s,t}$, linearity gives
\[
 \fieldnorm{(I+\mathcal T_\varepsilon Y_1)-(I+\mathcal T_\varepsilon Y_2)}
 =\fieldnorm{\mathcal T_\varepsilon(Y_1-Y_2)}
 \le\tfrac12\fieldnorm{Y_1-Y_2}.
\]
Thus $Y\mapsto I+\mathcal T_\varepsilon Y$ is a contraction of the Banach space
$\mathcal X_{s,t}$ into itself. Banach's fixed point theorem gives the
unique fixed point $Z^\varepsilon=I+\mathcal T_\varepsilon Z^\varepsilon$.
In particular, $\fieldnorm{Z^\varepsilon}\le2\fieldnorm{I}$ for every
$\varepsilon\ge0$, which gives the uniform second-moment bound.

\emph{(ii)} Write $Z\coloneqq Z^0$. Subtracting \eqref{eq:linear} at $\varepsilon=0$ from the same equation at
$\varepsilon>0$, we obtain
\[
 Z^\varepsilon-Z
 =\mathcal T_\varepsilon(Z^\varepsilon-Z)+(\mathcal T_\varepsilon-\mathcal T_0)Z .
\]
Fix $s<r\le t$ and take $Y=Z$ in the definition of $H_{q,x}$ in (i).
By \eqref{eq:smoothed-noise-integral} and the white-noise isometry, we obtain
\[
 \E\|[(\mathcal T_\varepsilon-\mathcal T_0)Z]_r\|_{2,\rho}^2
 =\int_s^r\!\int_\R\rho(x)\E\|P_\varepsilon H_{q,x}-H_{q,x}\|_2^2
          \,\d x\,\d q\longrightarrow0 \qquad(\varepsilon\downarrow0).
\]
Indeed, by strong continuity on $L^2(\R,\d y)$, the squared difference tends to zero
for almost every $(\omega,q,x)$. By contractivity it is bounded by
$4\|H_{q,x}\|_2^2$, which is integrable by \eqref{eq:smoothed-integrand-bound}.
It remains to control $\mathcal T_\varepsilon(Z^\varepsilon-Z)$.
For any sequence $\varepsilon_n\downarrow0$, set
$d(r)\coloneqq \limsup_n\E\|Z^{\varepsilon_n}_r-Z_r\|_{2,\rho}^2$.
Since the other term tends to zero in mean square, the $L^2$ triangle
inequality and \eqref{eq:linear-iteration} give
\[
 \begin{aligned}
 d(r)&\le K_\beta^2\limsup_{n\to\infty}\int_s^r
       \frac{\e^{c(r-q)/2}}{2\sqrt{\pi(r-q)}}\E\|Z_q^{\varepsilon_n}-Z_q\|_{2,\rho}^2\,\d q\\
 &\le K_\beta^2\int_s^r \frac{\e^{c(r-q)/2}}{2\sqrt{\pi(r-q)}}d(q)\,\d q .
 \end{aligned}
\]
The last inequality uses the integrable bound $C/\sqrt{r-q}$ from the uniform
second-moment estimate in (i).
Using the same exponential weight as in (i), we obtain
\[
 \sup_{s\le r\le t}\e^{-\upsilon(r-s)}d(r)
 \le\tfrac14\sup_{s\le r\le t}\e^{-\upsilon(r-s)}d(r).
\]
Since $d$ is nonnegative and bounded, we have $d(r)=0$ for every $r$.
The sequence $\varepsilon_n\downarrow0$ was arbitrary, which proves \eqref{eq:white}.
\end{proof}

\begin{lemma}\label{lem:smoothed-nonnegativity}
Under the assumptions of Proposition~\ref{prop:linear-comparison}, fix
$\varepsilon>0$ and let $Z^\varepsilon$ solve \eqref{eq:linear}.
If $\xi$ and $G$ are nonnegative, then $Z^\varepsilon_r\ge0$ for each $r\in[s,t]$.
\end{lemma}

\begin{proof}
To detect the negative part, set $z^-\coloneqq\max\{-z,0\}$ and
$\psi(z)\coloneqq (z^-)^3/(1+z^-)^2$.
Then $\psi\in C^2(\R)$, $\psi(z)=0$ exactly when $z\ge0$, and $-1\le\psi'\le0$.
Since $\psi''(z)=6z^-/(1+z^-)^4$, we have
\begin{equation}\label{eq:positivity-test}
 0\le z^2\psi''(z)=\frac{6\psi(z)}{(1+z^-)^2}\le6\psi(z).
\end{equation}

Let $K_q(x,y)\coloneqq\beta_q(x)Z^\varepsilon_q(x)p_\varepsilon(x-y)$.
By \eqref{eq:smoothed-noise-integral}, the stochastic integral in
\eqref{eq:linear} has integrand $P_{r-q}^xK_q$ against $W$, where $P_h^x$
acts on the first space variable. Since $\|p_\varepsilon\|_2^2=1/(2\sqrt{\pi\varepsilon})$,
Lemma~\ref{lem:linear-approximation}(i) gives
\[
 \E\int_s^t\!\int_\R\rho(x)\|K_q(x,\cdot)\|_2^2\,\d x\,\d q
 \le \frac{K_\beta^2}{2\sqrt{\pi\varepsilon}}\int_s^t\E\|Z^\varepsilon_q\|_{2,\rho}^2\,\d q<\infty.
\]
For $\delta>0$, set $\widetilde Z_q\coloneqq P_\delta Z^\varepsilon_q$ and
$K_q^\delta\coloneqq P_\delta^xK_q$. Applying $P_\delta$ to \eqref{eq:linear}
and using the semigroup property, we obtain
\[
 \widetilde Z_r(x)=(P_{\delta+r-s}\xi)(x)+\int_s^r(P_{\delta+r-q}G_q)(x)\,\d q
       +\int_s^r\!\int_\R (P_{\delta+r-q}^xK_q)(x,y)\,W(\d q,\d y).
\]
By the heat equation, we have
\[
 P_{\delta+r-q}=P_\delta+\tfrac12\int_q^r\partial_{xx}P_{\delta+\ell-q}\,\d\ell .
\]
Fix $x\in\R$. Substituting this identity and applying the stochastic
Fubini theorem \cite[Theorem~2.6]{Walsh},
we obtain a version of $(\widetilde Z_r(x))_{s\le r\le t}$ satisfying
\begin{equation}\label{eq:smoothed-semimartingale}
 \begin{aligned}
 \widetilde Z_r(x)&=(P_\delta\xi)(x)
   +\int_s^r\big(\tfrac12\partial_{xx}\widetilde Z_q(x)+(P_\delta G_q)(x)\big)\,\d q+N_r^x,\\
 N_r^x&\coloneqq\int_s^r\!\int_\R K_q^\delta(x,y)\,W(\d q,\d y).
 \end{aligned}
\end{equation}
For the spatial interchange, $f\mapsto(P_\delta f)(x)$ is bounded on
$L^2(\R,\rho(z)\d z)$ and Lemma~\ref{lem:heat-bounds} applies to $P_{r-q}^xK_q$.
For the time interchange, we use the square integrability of $K$ above and
$\sup_{\delta\le h\le\delta+t-s}\int_\R|\partial_{xx}p_h(x-z)|^2/\rho(z)\,\d z<\infty$.

The same pointwise bounds give an absolutely integrable drift and
$\E\int_s^t\|K_q^\delta(x,\cdot)\|_2^2\,\d q<\infty$.
Thus $N^x$ is a real-valued continuous square-integrable martingale with
$[N^x]_r=\int_s^r\|K_q^\delta(x,\cdot)\|_2^2\,\d q$
(see \cite[Theorem~2.5, (2.10)]{Walsh}). Applying the one-dimensional It\^o
formula to the continuous semimartingale in \eqref{eq:smoothed-semimartingale}, we obtain
\[
 \begin{aligned}
 \psi(\widetilde Z_r(x))-\psi(\widetilde Z_s(x))
 &=\int_s^r\psi'(\widetilde Z_q(x))
    \big(\tfrac12\partial_{xx}\widetilde Z_q(x)+(P_\delta G_q)(x)\big)\,\d q
    +\int_s^r\psi'(\widetilde Z_q(x))\,\d N_q^x\\
 &\quad+\tfrac12\int_s^r\psi''(\widetilde Z_q(x))\|K_q^\delta(x,\cdot)\|_2^2\,\d q .
 \end{aligned}
\]
We now take expectations and integrate against $\rho(x)\d x$. The integral
against $N^x$ has mean zero since $|\psi'|\le1$. We have
$\psi(P_\delta\xi)=0$ and $\psi'(\widetilde Z_q)P_\delta G_q\le0$.
Using integration by parts and \eqref{eq:rho-generator-bound}, we obtain
\[
 \tfrac12\langle\rho,\psi'(\widetilde Z_q)\partial_{xx}\widetilde Z_q\rangle
 =\tfrac12\langle\rho'',\psi(\widetilde Z_q)\rangle
  -\tfrac12\langle\rho,\psi''(\widetilde Z_q)|\partial_x \widetilde Z_q|^2\rangle
 \le c\langle\rho,\psi(\widetilde Z_q)\rangle .
\]
Here we may remove compact cutoffs since $|\rho'|+|\rho''|\le C_\lambda\rho$
and $|\partial_x^j p_\delta|\le C_{\delta,j}p_{2\delta}$ for $j=1,2$.
It follows that
\[
 \E\langle\rho,\psi(\widetilde Z_r)\rangle
 \le c\int_s^r\E\langle\rho,\psi(\widetilde Z_q)\rangle\,\d q
  +\tfrac12\E\int_s^r\!\int_\R\rho(x)\psi''(\widetilde Z_q(x))
                         \|K_q^\delta(x,\cdot)\|_2^2\,\d x\,\d q .
\]

By density and Lemma~\ref{lem:heat-bounds}, $P_\delta$ converges strongly
to the identity in $L^2(\R,\rho(x)\d x)$, also for $L^2(\R,\d y)$-valued functions.
Thus $\widetilde Z_q\to Z^\varepsilon_q$ in weighted mean square at each $q$ and after
integration in $q$. Likewise $K^\delta\to K$ in weighted mean square
after integration in $q$, so $\|K_q^\delta(x,\cdot)\|_2^2\to\|K_q(x,\cdot)\|_2^2$
in $L^1(\d\mathbb P\,\d q\,\rho(x)\d x)$. Since $\psi''$ is bounded and
continuous, $\psi''(\widetilde Z)\to\psi''(Z^\varepsilon)$ in measure, which justifies
passing to the limit in the It\^o correction. Set
$m(r)\coloneqq\E\langle\rho,\psi(Z^\varepsilon_r)\rangle$. Letting $\delta\downarrow0$
and using $\|K_q(x,\cdot)\|_2^2=|\beta_q(x)Z^\varepsilon_q(x)|^2/(2\sqrt{\pi\varepsilon})$
and \eqref{eq:positivity-test}, we obtain
\[
 m(r)\le\left(c+\frac{3K_\beta^2}{2\sqrt{\pi\varepsilon}}\right)\int_s^r m(q)\,\d q .
\]
Gronwall's inequality gives $m(r)=0$. Since $\rho>0$ and $\psi(z)>0$ for
$z<0$, we obtain $Z^\varepsilon_r\ge0$.
\end{proof}

\begin{proof}[Proof of Proposition~\ref{prop:linear-comparison}]
Since $M^0=W$, equation \eqref{eq:linear} at $\varepsilon=0$ coincides with
\eqref{eq:linear-white}. Lemma~\ref{lem:linear-approximation}(i) therefore gives
the asserted existence and uniqueness, with $Z=Z^0$.

\emph{(i)} Fix $r\in[s,t]$. For $\xi,G\ge0$, Lemma~\ref{lem:smoothed-nonnegativity}
gives $Z^\varepsilon_r\ge0$ for every $\varepsilon>0$.
Recall that $Z_r^-$ denotes the negative part of $Z_r$.
By \eqref{eq:white}, we have
\[
 \E\|Z_r^-\|_{2,\rho}^2
 \le \E\|Z_r-Z^\varepsilon_r\|_{2,\rho}^2\longrightarrow0 .
\]
Hence $Z_r\ge0$.

\emph{(ii)} For signed data, apply the existence assertion just proved to
\eqref{eq:linear-white} with initial value $|\xi|$, drift $|G|$ and the
same $\beta$ and $W$. Denote the resulting solution by $\overline Z$.
By linearity, $\overline Z\pm Z$ solve the equations with nonnegative data
$|\xi|\pm\xi$ and $|G|\pm G$. Applying (i) to both gives $|Z_r|\le \overline Z_r$.
Pairing the mild equation for $\overline Z_r$ with $\rho$ and taking expectations, we obtain
\[
 \begin{aligned}
 \E\|Z_r\|_{1,\rho}
 &\le\E\langle\rho,\overline Z_r\rangle\\
 &=\E\langle P_{r-s}\rho,|\xi|\rangle
   +\int_s^r\E\langle P_{r-q}\rho,|G_q|\rangle\,\d q .
 \end{aligned}
\]
By Lemma~\ref{lem:stochastic-convolution}, the stochastic convolution is centered
and square integrable in $L^2(\R,\rho(x)\d x)$. Hence its pairing with $\rho$ has
mean zero, since $|\langle\rho,f\rangle|\le\|\rho\|_1^{1/2}\|f\|_{2,\rho}$.
Applying the first bound in Lemma~\ref{lem:heat-bounds} gives \eqref{eq:source-mass}.
\end{proof}

\section{Stability and pathwise uniqueness}\label{sec:block}

{
Let $u$ be a solution in the sense of Definition~\ref{def:solution}(i)
with initial function $u_0$ and driving noise $W$. Then, for each fixed $r>0$
and $x\in\R$, the equivalence of the weak and mild formulations
\cite[Theorem~2.1]{Shiga} gives the almost sure identity
\[
 \begin{aligned}
 u_r(x)={}&(P_r u_0)(x)
 +\int_0^r\!\int_\R p_{r-\tau}(x-y)b(\tau,y,u_\tau(y))\,\d y\,\d \tau\\
 &+\int_0^r\!\int_\R p_{r-\tau}(x-y)\sigma(\tau,y,u_\tau(y))\,W(\d\tau,\d y).
 \end{aligned}
\]
For a predictable field $B$ for which the integral is defined, write
\[
 \mathcal W_{s,r}(B)(x)\coloneqq
 \int_s^r\!\int_\R p_{r-q}(x-y)B_q(y)\,W(\d q,\d y),\qquad 0\le s\le r .
\]
On $[s,t]$, let $Z$ be the auxiliary process with initial value $\xi$,
drift $G$ and noise coefficient $\beta$ from Proposition~\ref{prop:linear-comparison}.
We can now write its mild equation \eqref{eq:linear-white} as
\[
 Z_r=P_{r-s}\xi+\int_s^r P_{r-q}G_q\,\d q
             +\mathcal W_{s,r}(\beta Z),\qquad s\le r\le t.
\]}

\begin{proposition}\label{prop:stability}
Suppose that \condH{} holds. Let $u,v$ be solutions
with initial functions $u_0,v_0\in C_{\rm tem}$ on the same filtered
space and driven by the same white noise. Fix $T<\infty$ and take the constants
$L,a$ in \condH{} on $[0,T]$. Then, for every $\lambda>8a$ and the
corresponding weight $\rho$, we have, for $0\le t\le T$,
\begin{equation}\label{eq:stability}
 \E\|u_t-v_t\|_{1,\rho}
 \le \e^{(\lambda^2/2+L)t}\|u_0-v_0\|_{1,\rho}.
\end{equation}
\end{proposition}

\begin{proof}[Proof of Proposition~\ref{prop:stability}]
\emph{Step 1.} We localize the two solutions and establish the
integrability bounds needed for the comparison argument.

Fix $0<t\le T$ and $\lambda>8a$.
We may take $K\ge2C$ and $m\ge1$ in \condH{}.
Choose $\kappa>0$ sufficiently small that
\begin{equation}\label{eq:kappa-choice}
 8(a+m\kappa)\mathbin{+}2\kappa<\lambda .
\end{equation}
For the two solutions $u,v$ in the proposition, define the coefficient
differences by
\[
 A_r(x)=\sigma(r,x,u_r(x))-\sigma(r,x,v_r(x)),\qquad
 G_r(x)=b(r,x,u_r(x))-b(r,x,v_r(x)).
\]
For each integer $n>|u_0|_\kappa\vee|v_0|_\kappa$, define the stopping
time and the stopped coefficients by
\[
 \tau_n=T\wedge\inf\{r\ge0:|u_r|_\kappa\vee|v_r|_\kappa\ge n\},
 \qquad(A^n_r,G^n_r)=\ind_{\{r\le\tau_n\}}(A_r,G_r).
\]
For $w=u,v$, define
\[
 w^n_r\coloneqq
 \begin{cases}
  w_r,&0\le r\le\tau_n,\\
  P_{r-\tau_n}w_{\tau_n},&\tau_n<r\le T.
 \end{cases}
\]
Using the weak equation up to $\tau_n$ and the heat equation afterwards,
we obtain the mild form
\[
 \begin{aligned}
 w^n_r(x)
 &=(P_rw_0)(x)
   +\int_0^{r\wedge\tau_n}\!\int_\R
       p_{r-q}(x-y)b(q,y,w^n_q(y))\,\d y\,\d q\\
 &\quad+\int_0^{r\wedge\tau_n}\!\int_\R
       p_{r-q}(x-y)\sigma(q,y,w^n_q(y))\,W(\d q,\d y),
       \qquad 0\le r\le T.
 \end{aligned}
\]
Set $D^n\coloneqq u^n-v^n$. Subtracting the two mild equations and using
$u^n=u$, $v^n=v$ up to $\tau_n$, we obtain
\begin{equation}\label{eq:stopped}
 D^n_r=P_r(u_0-v_0)+\int_0^rP_{r-\tau}G^n_\tau\d \tau
                     +\mathcal W_{0,r}(A^n),\qquad 0\le r\le T.
\end{equation}
By construction, almost surely for all
$(r,x)\in[0,T]\times\R$, we have
\begin{equation}\label{eq:stop-equal}
 \ind_{\{r\le\tau_n\}}D^n_r(x)
 =\ind_{\{r\le\tau_n\}}\bigl(u_r(x)-v_r(x)\bigr).
\end{equation}
On $\{r\le\tau_n\}$, we have
$|u_r(x)|,\,|v_r(x)|\le n\e^{\kappa|x|}$, so
$1+|u_r(x)|+|v_r(x)|\le(1+2n)\e^{\kappa|x|}$.
Set
\begin{equation}\label{eq:Jn-definition}
 J_n(x)\coloneqq K(1+2n)^m\e^{(a+m\kappa)|x|}.
\end{equation}
Using \condH{}(i), $K\ge2C$ and $m\ge1$, we obtain on this event
\begin{equation}\label{eq:stopped-linear-growth}
 |A^n_r(x)|+|G^n_r(x)|
 \le 2C(1+n\e^{\kappa|x|})\le J_n(x).
\end{equation}
On $\{r\le\tau_n,\ |D^n_r(x)|\le1\}$, applying \condH{}(iii) and
\eqref{eq:stop-equal}, we obtain
\begin{equation}\label{eq:localized-holder}
 \begin{aligned}
 |A^n_r(x)|
 &\le K\e^{a|x|}(1+|u_r(x)|+|v_r(x)|)^m|D^n_r(x)|^\gamma\\
 &\le J_n(x)|D^n_r(x)|^\gamma.
 \end{aligned}
\end{equation}
Combining these estimates and using $A^n=G^n=0$ after $\tau_n$, we obtain
almost surely for all $(r,x)\in[0,T]\times\R$
\begin{equation}\label{eq:stopped-bounds}
 |A^n_r(x)|\le J_n(x)\bigl(1\wedge|D^n_r(x)|^\gamma\bigr),\qquad
 |G^n_r(x)|\le L|D^n_r(x)|,
\end{equation}
where the second inequality follows from \condH{}(ii).

\Needspace{9\baselineskip}
{Using \eqref{eq:Jn-definition} and $\sqrt{1+x^2}\ge|x|$, we obtain
\[
 \begin{aligned}
 \rho(x)J_n(x)^8\e^{2\kappa|x|}
 &=\e^{-\lambda\sqrt{1+x^2}}
   \bigl[K(1+2n)^m\e^{(a+m\kappa)|x|}\bigr]^8\e^{2\kappa|x|}\\
 &\le K^8(1+2n)^{8m}\e^{-[\lambda-8(a+m\kappa)-2\kappa]|x|}.
 \end{aligned}
\]
By \eqref{eq:kappa-choice}, $\lambda-8(a+m\kappa)-2\kappa>0$.
Integrating this bound over $\R$ therefore gives
\begin{equation}\label{eq:Jn-weight}
 \int_\R\rho(x)J_n(x)^8\,\d x
 \le \int_\R\rho(x)J_n(x)^8\e^{2\kappa|x|}\,\d x
 \le \frac{2K^8(1+2n)^{8m}}{\lambda-8(a+m\kappa)-2\kappa}<\infty .
\end{equation}}
As $\rho\in L^1(\R,\d x)$, we also have $\|J_n\|_{2,\rho}<\infty$.
For $0\le r\le T$, applying Lemma~\ref{lem:heat-bounds} to the defining
formula for $D^n$ gives
\[
 \begin{aligned}
 \|D^n_r\|_{2,\rho}^2
 &=\bigl\|P_{r-(r\wedge\tau_n)}
       (u_{r\wedge\tau_n}-v_{r\wedge\tau_n})\bigr\|_{2,\rho}^2\\
 &\le \e^{cT}\|u_{r\wedge\tau_n}-v_{r\wedge\tau_n}\|_{2,\rho}^2\\
 &\le 4n^2\e^{cT}\int_\R\rho(x)\e^{2\kappa|x|}\,\d x<\infty.
 \end{aligned}
\]
Here $2\kappa<\lambda$ by \eqref{eq:kappa-choice}, and the final bound is deterministic and uniform in $r$.
Using $|A^n|+|G^n|\le J_n$, we also obtain
\[
 \begin{aligned}
 \|A^n_r\|_{2,\rho}^2+\|G^n_r\|_{2,\rho}^2
 &=\int_\R\rho(x)\bigl(|A^n_r(x)|^2+|G^n_r(x)|^2\bigr)\,\d x\\
 &\le\|J_n\|_{2,\rho}^2<\infty.
 \end{aligned}
\]
Taking expectations and then the supremum over $r\le T$ in these bounds,
we conclude that
\begin{equation}\label{eq:localized-moments}
 \sup_{r\le T}\E\big[
 \|A^n_r\|_{2,\rho}^2+\|G^n_r\|_{2,\rho}^2
                         +\|D^n_r\|_{2,\rho}^2\big]<\infty.
\end{equation}
{
\smallskip\noindent\emph{Step 2.} We construct the auxiliary processes and
reduce the desired estimate to controlling the sum of the comparison errors.

Fix $n$. Put $h\coloneqq t/N$ and $t_j\coloneqq jh$.
To construct the auxiliary processes with a bounded coefficient, define
\[
 \beta\coloneqq
 \begin{cases}
  A^n/D^n,&D^n\ne0,\\
  0,&D^n=0,
 \end{cases}
 \qquad
 \beta^{(h)}\coloneqq
 \begin{cases}
  -\tfrac12h^{-1/4},&\beta<-\tfrac12h^{-1/4},\\
  \beta,&|\beta|\le\tfrac12h^{-1/4},\\
  \tfrac12h^{-1/4},&\beta>\tfrac12h^{-1/4}.
 \end{cases}
\]
For $j=0,\ldots,N-1$, consider the equation
\[
 Z^j_r=P_{r-t_j}D^n_{t_j}
       +\int_{t_j}^rP_{r-q}G^n_q\,\d q
       +\mathcal W_{t_j,r}(\beta^{(h)} Z^j),\qquad t_j\le r\le t_{j+1}.
\]
This is \eqref{eq:linear-white} with initial value $D^n_{t_j}$, drift $G^n$
and coefficient $\beta^{(h)}$. By \eqref{eq:localized-moments}, we have
\[
 \E\|D^n_{t_j}\|_{2,\rho}^2<\infty,
 \qquad
 \int_{t_j}^{t_{j+1}}\E\|G^n_q\|_{2,\rho}^2\,\d q
 \le h\sup_{q\le T}\E\|G^n_q\|_{2,\rho}^2<\infty.
\]
\begin{samepage}
Using these bounds and $|\beta^{(h)}|\le\tfrac12h^{-1/4}$, we apply
Proposition~\ref{prop:linear-comparison} to obtain a unique predictable
solution $Z^j$ satisfying
\begin{equation}\label{eq:comparison-moments}
 \sup_{t_j\le r\le t_{j+1}}\E\|Z^j_r\|_{2,\rho}^2<\infty,
\end{equation}
and
\begin{equation}\label{eq:comparison-l1}
 \E\|Z^j_{t_{j+1}}\|_{1,\rho}
 \le \e^{ch}\E\|D^n_{t_j}\|_{1,\rho}
 +\int_{t_j}^{t_{j+1}}\e^{c(t_{j+1}-q)}\E\|G^n_q\|_{1,\rho}\,\d q.
\end{equation}
\end{samepage}
To compare $D^n$ with $Z^j$, set $R_h\coloneqq A^n-\beta^{(h)}D^n$ and
$X_r\coloneqq D^n_r-Z^j_r$ for $t_j\le r\le t_{j+1}$.
Restarting \eqref{eq:stopped} at $t_j$ and subtracting the mild equation
for $Z^j$, we cancel the initial and drift terms and obtain
\begin{equation}\label{eq:error}
 X_r=\mathcal W_{t_j,r}(\beta^{(h)} X+R_h),\qquad t_j\le r\le t_{j+1}.
\end{equation}
To transfer the $L^1$ estimate to $D^n$, set
\[
 \varepsilon_j\coloneqq\E\|X_{t_{j+1}}\|_{1,\rho}.
\]
Applying \eqref{eq:comparison-l1} and the second bound in
\eqref{eq:stopped-bounds}, we obtain
\[
 \begin{aligned}
 \E\|D^n_{t_{j+1}}\|_{1,\rho}
 &\le\E\|Z^j_{t_{j+1}}\|_{1,\rho}+\varepsilon_j\\
 &\le\e^{ch}\E\|D^n_{t_j}\|_{1,\rho}
 +L\int_{t_j}^{t_{j+1}}\e^{c(t_{j+1}-q)}\E\|D^n_q\|_{1,\rho}\,\d q
 +\varepsilon_j.
 \end{aligned}
\]
Multiplying by $\e^{-ct_{j+1}}$ and summing over $j$, we obtain
\begin{equation}\label{eq:comparison-recursion}
 \e^{-ct}\E\|D^n_t\|_{1,\rho}
 \le\|u_0-v_0\|_{1,\rho}
       +L\int_0^t\e^{-cq}\E\|D^n_q\|_{1,\rho}\,\d q
       +\sum_{j=0}^{N-1}\e^{-ct_{j+1}}\varepsilon_j.
\end{equation}
For fixed $n$, the process $D^n$ and the constants $c,L$ do not depend on $N$.
Thus, to obtain a Gronwall inequality, it suffices to prove
\begin{equation}\label{eq:comparison-error-limit}
 \sum_{j=0}^{N-1}\varepsilon_j\longrightarrow0
 \qquad\text{as }N\to\infty.
\end{equation}

\smallskip\noindent\emph{Step 3.} We next show that
\eqref{eq:comparison-error-limit} holds and complete the proof.

On the truncation event $\{|\beta_q(x)|>\tfrac12 h^{-1/4}\}$,
we have $D^n_q(x)\ne0$ by the definition of $\beta_q(x)$.
Since $0\le4\gamma-3\le1$, \eqref{eq:stopped-linear-growth} and
\eqref{eq:stopped-bounds} give
\begin{equation}\label{eq:truncation-coefficient-bound}
 |\beta|^{4\gamma}|A^n|
 =\left(\frac{|A^n|}{|D^n|^\gamma}\right)^4 |A^n|^{4\gamma-3}
 \le J_n(x)^4 |A^n|^{4\gamma-3}
 \le C_n J_n(x)^4\e^{\kappa|x|}.
\end{equation}
Using the definitions of $\beta^{(h)}$ and $R_h$
and applying \eqref{eq:truncation-coefficient-bound} in the last line, we obtain
\begin{equation}\label{eq:critical-clip}
 \begin{aligned}
 h^{-\gamma}|R_h|
 &=h^{-\gamma}|D^n|\,|\beta-\beta^{(h)}|\\
 &=h^{-\gamma}|D^n|\bigl(|\beta|-\tfrac12h^{-1/4}\bigr)^+\\
 &\le h^{-\gamma}|A^n|\ind_{\{|\beta|>\frac12h^{-1/4}\}}\\
 &\le 2^{4\gamma}|\beta|^{4\gamma}|A^n|\ind_{\{|\beta|>\frac12h^{-1/4}\}}\\
 &\le C_n J_n(x)^4\e^{\kappa|x|}\ind_{\{|\beta|>\frac12h^{-1/4}\}}.
 \end{aligned}
\end{equation}

To bound $\varepsilon_j$ on $I_j=[t_j,t_{j+1}]$, we first square
\eqref{eq:critical-clip} to obtain
\[
 h^{-2\gamma}\E\|R_h(q)\|_{2,\rho}^2
 \le C_n\,\E\int_\R\rho(x)J_n(x)^8\e^{2\kappa|x|}
       \ind_{\{|\beta_q(x)|>\frac12h^{-1/4}\}}\,\d x .
\]
The indicator tends to zero pointwise as $h\downarrow0$.
By \eqref{eq:Jn-weight} and dominated convergence,
\begin{equation}\label{eq:residual-moment-rates}
 \sup_{0\le q\le t}\E\|R_h(q)\|_{2,\rho}^2=O(h^{2\gamma}),\qquad
 \int_0^t\E\|R_h(q)\|_{2,\rho}^2\,\d q=o(h^{2\gamma}).
\end{equation}

\enlargethispage{\baselineskip}
Take $h$ small enough that $\e^{ch/2}\le\sqrt\pi$.
Applying \eqref{eq:isometry} to \eqref{eq:error} gives, for $r\in I_j$,
\[
 \begin{aligned}
 \E\|X_r\|_{2,\rho}^2
 &\le \int_{t_j}^r \frac{\e^{c(r-q)/2}}{2\sqrt{\pi(r-q)}}\,
       \E\|\beta_q^{(h)}X_q+R_h(q)\|_{2,\rho}^2\,\d q\\
 &\le \int_{t_j}^r
       \frac{\E\|\beta_q^{(h)}X_q\|_{2,\rho}^2+\E\|R_h(q)\|_{2,\rho}^2}
            {\sqrt{r-q}}\,\d q\\
 &\le \underbrace{\bigl(\tfrac12h^{-1/4}\bigr)^2
       \int_{t_j}^r\frac{\d q}{\sqrt{r-q}}}
       _{=\,\frac12\sqrt{(r-t_j)/h}\,\le\,1/2}\,
       \sup_{s\in I_j}\E\|X_s\|_{2,\rho}^2
     +\int_{t_j}^r\frac{\E\|R_h(q)\|_{2,\rho}^2}{\sqrt{r-q}}\,\d q ,
 \end{aligned}
\]
where the last inequality uses $|\beta_q^{(h)}|\le\tfrac12 h^{-1/4}$.

The supremum is finite by \eqref{eq:localized-moments} and
\eqref{eq:comparison-moments}. Squaring the last integral and using symmetry gives
\[
 \begin{aligned}
 &\left(\int_{t_j}^r\frac{\E\|R_h(q)\|_{2,\rho}^2}{\sqrt{r-q}}\,\d q\right)^2\\
 &\quad=2\int_{t_j}^r\frac{\E\|R_h(q)\|_{2,\rho}^2}{\sqrt{r-q}}\,
       \underbrace{\int_q^r\frac{\E\|R_h(v)\|_{2,\rho}^2}{\sqrt{r-v}}\,\d v}
          _{\le\,2\sqrt{r-q}\,\sup_{0\le s\le t}\E\|R_h(s)\|_{2,\rho}^2}\,\d q\\
 &\quad\le 4\left(\sup_{0\le s\le t}\E\|R_h(s)\|_{2,\rho}^2\right)
           \int_{I_j}\E\|R_h(q)\|_{2,\rho}^2\,\d q .
 \end{aligned}
\]
Taking the supremum over $r\in I_j$ in the estimate for $\E\|X_r\|_{2,\rho}^2$
and moving the first term to the left gives
\begin{equation}\label{eq:comparison-second-moment}
 \frac12\sup_{r\in I_j}\E\|X_r\|_{2,\rho}^2
 \le 2\left(\sup_{0\le s\le t}\E\|R_h(s)\|_{2,\rho}^2
              \int_{I_j}\E\|R_h(q)\|_{2,\rho}^2\,\d q\right)^{1/2}.
\end{equation}
By Cauchy--Schwarz,
\[
 \varepsilon_j^2
 =\bigl(\E\|X_{t_{j+1}}\|_{1,\rho}\bigr)^2
 \le\|\rho\|_1\,\E\|X_{t_{j+1}}\|_{2,\rho}^2 .
\]
Applying \eqref{eq:comparison-second-moment} gives
\[
 \varepsilon_j
 \le 2\|\rho\|_1^{1/2}
       \left[\left(\sup_{0\le s\le t}\E\|R_h(s)\|_{2,\rho}^2\right)
             \int_{I_j}\E\|R_h(q)\|_{2,\rho}^2\,\d q\right]^{1/4}.
\]
Summing over the $N=t/h$ intervals and using the concavity of $z\mapsto z^{1/4}$, we obtain
\begin{equation}\label{eq:sum}
 \begin{aligned}
 \sum_{j=0}^{N-1}\varepsilon_j
 &\le 2\|\rho\|_1^{1/2}
       \left(\sup_{0\le s\le t}\E\|R_h(s)\|_{2,\rho}^2\right)^{1/4}
       \sum_{j=0}^{N-1}\left(\int_{I_j}\E\|R_h(q)\|_{2,\rho}^2\,\d q\right)^{1/4}\\
 &\le 2\|\rho\|_1^{1/2}
       \left(\sup_{0\le s\le t}\E\|R_h(s)\|_{2,\rho}^2\right)^{1/4}
       N\left(\frac1N\sum_{j=0}^{N-1}\int_{I_j}\E\|R_h(q)\|_{2,\rho}^2\,\d q\right)^{1/4}\\
 &= 2\|\rho\|_1^{1/2}
       \left(\frac{t}{h}\right)^{3/4}
       \Biggl[\underbrace{\left(\sup_{0\le s\le t}\E\|R_h(s)\|_{2,\rho}^2\right)}_{O(h^{2\gamma})}
             \,\underbrace{\int_0^t\E\|R_h(q)\|_{2,\rho}^2\,\d q}_{o(h^{2\gamma})}\Biggr]^{1/4}
       =o\bigl(h^{\gamma-3/4}\bigr)
 \end{aligned}
\end{equation}
as $h\downarrow0$, by \eqref{eq:residual-moment-rates}.}

Applying \eqref{eq:sum} in \eqref{eq:comparison-recursion}
and letting $N\to\infty$ with $n$ and $t$ fixed, we obtain
\[
 \e^{-ct}\E\|D^n_t\|_{1,\rho}
 \le\|u_0-v_0\|_{1,\rho}
 +L\int_0^t\e^{-cq}\E\|D^n_q\|_{1,\rho}\,\d q,
 \qquad 0\le t\le T.
\]
By \eqref{eq:localized-moments}, $r\mapsto\E\|D^n_r\|_{1,\rho}$
is bounded on $[0,T]$, so Gronwall's inequality gives
\[
 \E\|D^n_t\|_{1,\rho}\le\e^{(c+L)t}\|u_0-v_0\|_{1,\rho}.
\]
{Almost surely the paths are bounded in $|\cdot|_\kappa$ on $[0,T]$,
so eventually $\tau_n=T$ and $D^n=u-v$ on $[0,T]$ by \eqref{eq:stop-equal}.
Fatou's lemma gives
\[
 \E\|u_t-v_t\|_{1,\rho}
 \le\liminf_{n\to\infty}\E\|D^n_t\|_{1,\rho}
 \le\e^{(c+L)t}\|u_0-v_0\|_{1,\rho}.
\]}
If $u_0=v_0$, spatial
continuity gives equality for every $x$ at each rational time almost
surely, and time continuity gives indistinguishability.
\end{proof}

\section{Existence and uniqueness in law}\label{sec:existence}

\begingroup

First, we regularize the diffusion coefficient and construct
approximate solutions by Picard iteration.
Lemma~\ref{lem:ex-approx} gives the convergence
of the iteration and uniform moment bounds.
Second, Lemma~\ref{lem:ex-weak} proves weak existence.
Finally, we establish Theorem~\ref{thm:main} by combining weak existence from
Lemma~\ref{lem:ex-weak} with pathwise uniqueness
from Proposition~\ref{prop:stability}.

Let $\pi_n z=(-n)\vee(z\wedge n)$.
Choose a smooth probability density $\varphi$ supported on $[-1,1]$.
Set $\varphi_n(v):=n\varphi(nv)$ and put
\[
 \sigma_n(r,x,z):=\int_\R\varphi_n(v)\sigma(r,x,\pi_n(z-v))\d v .
\]
Fix $u_0\in C_{\rm tem}$ and $T>0$.
By \condH{}(i) and (iii), we have
\begin{equation}\label{eq:ex-mollifier}
 |\sigma_n(r,x,z)|\le C(2+|z|),\qquad
 |\partial_z\sigma_n(r,x,z)|\le C(1+n)\|\varphi_n'\|_1
\end{equation}
and
\begin{align*}
 \sup_{\substack{r\le T,~|x|,|z|\le R}}|\sigma_n(r,x,z)-\sigma(r,x,z)|
 \le K\e^{aR}(2+2R)^m n^{-\gamma},\quad n>R+1.
\end{align*}
Let $W(\d s,\d u)$ be a space-time white noise with intensity $\d s \d u$.

\begin{lemma}\label{lem:ex-approx}
For each $n\ge1$, equation \eqref{eq:spde} with $\sigma$ replaced
by $\sigma_n$ admits a predictable mild solution $u^n$ with initial
condition $u_0$, driven by the given white noise $W$.
Moreover, for each $T>0$, $\mu>0$ and $p>8$,
\begin{equation}\label{eq:ex-moment}
 \sup_{n\ge1,\,0\le t\le T,\,x\in\R}
 \e^{-2\mu|x|}\|u^n_t(x)\|_{L^p(\Omega)}^2<\infty.
\end{equation}
\end{lemma}

\begin{proof}
Fix $p>8$ and $\mu>0$. For each $n$, set $V_n^0(t,x)=P_tu_0(x)$,
with $P_0u_0=u_0$. For $j\ge0$, set $V_n^{j+1}(0,x)=u_0(x)$.
For $0<t\le T$, define
\begin{equation}\label{eq:ex-picard-recursion}
\begin{aligned}
 V_n^{j+1}(t,x)
 &:= P_tu_0(x)
 +\int_0^t\int_\R p_{t-v}(x-y)b(v,y,V_n^j(v,y))\,\d y\,\d v\\
 &\quad+\int_0^t\int_\R p_{t-v}(x-y)\sigma_n(v,y,V_n^j(v,y))\,W(\d v,\d y).
\end{aligned}
\end{equation}
Let
\[
 D_{n,j}(t)=\sup_x\e^{-2\mu|x|}
             \|V_n^{j+1}(t,x)-V_n^j(t,x)\|_{L^p(\Omega)}^2.
\]
Let
$L_n=C(1+n)\|\varphi_n'\|_1$. For $j\ge1$, define the coefficient differences
\begin{align*}
 \Delta b_v(y)&=b(v,y,V_n^j(v,y))-b(v,y,V_n^{j-1}(v,y)),\\
 \Delta\sigma_v(y)&=\sigma_n(v,y,V_n^j(v,y))-\sigma_n(v,y,V_n^{j-1}(v,y)).
\end{align*}
It follows from \condH{}(ii) and \eqref{eq:ex-mollifier} that
\begin{align*}
 \|\Delta b_v(y)\|_{L^p(\Omega)}
 &\le L\|V_n^j(v,y)-V_n^{j-1}(v,y)\|_{L^p(\Omega)}
 \le L\e^{\mu|y|}D_{n,j-1}(v)^{1/2},\\
 \|\Delta\sigma_v(y)\|_{L^p(\Omega)}
 &\le L_n\e^{\mu|y|}D_{n,j-1}(v)^{1/2}.
\end{align*}
By \eqref{eq:ex-picard-recursion},
\begin{equation}\label{eq:ex-picard-split}
V_n^{j+1}(t,x)-V_n^j(t,x)=I_t(x)+J_t(x),
\end{equation}
where
\begin{align*}
 I_t(x)&=\int_0^t\int_\R p_{t-v}(x-y)\Delta b_v(y)\d y\,\d v,\\
 J_t(x)&=\int_0^t\int_\R p_{t-v}(x-y)\Delta\sigma_v(y)W(\d v,\d y).
\end{align*}
Observe that there is a constant $\tilde{C}_T>0$ such that
\[
\|\e^{\mu|\cdot|}p_s\|_1
+s^{1/2}\|\e^{\mu|\cdot|}p_s\|_2^2\le
\tilde{C}_T,\qquad   0<s\le T.
\]
Then
Minkowski's inequality and Cauchy-Schwarz's inequality give
\begin{equation}\label{eq:ex-drift-diff}
 \e^{-\mu|x|}\|I_t(x)\|_{L^p(\Omega)}
 \le L\int_0^t\|\e^{\mu|\cdot|}p_{t-v}\|_1 D_{n,j-1}(v)^{1/2}\d v
\le L\tilde{C}_TT^{1/2}\Big[\int_0^t D_{n,j-1}(v)\d v\Big]^{1/2}.
\end{equation}
Applying the Burkholder--Davis--Gundy inequality
and then Minkowski's inequality one obtains
\begin{equation}\label{eq:ex-noise-diff}
\begin{aligned}
 \e^{-2\mu|x|}\|J_t(x)\|_{L^p(\Omega)}^2
&\le C_p
 \e^{-2\mu|x|}\left\|\int_0^t\int_{\mathbb R}
p_{t-v}(x-y)^2|\Delta\sigma_v(y)|^2\,\mathrm dy\,\mathrm dv
\right\|_{L^{p/2}(\Omega)} \\
 &\le C_p\e^{-2\mu|x|}\int_0^t\int_\R
 p_{t-v}(x-y)^2\|\Delta\sigma_v(y)\|_{L^p(\Omega)}^2\d y\,\d v\\
 &\le C_p L_n^2\int_0^t\|\e^{\mu|\cdot|}p_{t-v}\|_2^2
                 D_{n,j-1}(v)\d v\\
 &\le C_p L_n^2 \tilde{C}_T\int_0^t(t-v)^{-1/2}D_{n,j-1}(v)\d v
 \end{aligned}
\end{equation}
for some constant $C_p>0$.
Combining \eqref{eq:ex-drift-diff} and \eqref{eq:ex-noise-diff} gives, for some $C_n$ independent of $j$ and $t$,
\begin{equation}\label{eq:ex-picard-increments}
 D_{n,j}(t)\le
 C_n\int_0^t\bigl(1+(t-v)^{-1/2}\bigr)D_{n,j-1}(v)\d v,
 \qquad j\ge1,\quad 0\le t\le T.
\end{equation}

Let $n\ge1$ be fixed. Define
\[
 \|X\|_{\beta,T}
 :=\sup_{0\le t\le T,\,x\in\R}
   \e^{-\beta t/2-\mu|x|}\|X(t,x)\|_{L^p(\Omega)}.
\]
Multiplying \eqref{eq:ex-picard-increments} by $\e^{-\beta t}$ and taking the supremum gives
\[
 \|V_n^{j+1}-V_n^j\|_{\beta,T}^2
 \le q\,\|V_n^j-V_n^{j-1}\|_{\beta,T}^2,\qquad
 q=C_n\int_0^\infty\e^{-\beta s}(1+s^{-1/2})\,\d s
   =C_n(\beta^{-1}+\sqrt{\pi/\beta}).
\]
Choose $\beta$ so that $q<1$. The linear growth bounds and
$|P_vu_0(y)|\le\tilde C_T|u_0|_\mu\e^{\mu|y|}$ give
$\|V_n^1-V_n^0\|_{\beta,T}<\infty$ by the same convolution estimates.
Thus the successive differences are summable in $\|\cdot\|_{\beta,T}$
and absolutely summable almost surely at each $(t,x)$. Their pointwise sum,
added to $V_n^0$ and defined as zero where the series diverges, gives a
predictable field $u^n$ with $\|V_n^j-u^n\|_{\beta,T}\to0$.
Applying \eqref{eq:ex-drift-diff} and \eqref{eq:ex-noise-diff} to the
coefficient differences between $V_n^j$ and $u^n$ shows that both
convolutions in \eqref{eq:ex-picard-recursion} converge in $L^p(\Omega)$.
Hence $u^n$ satisfies the mild equation with coefficients $(b,\sigma_n)$.
The recursion uses the same noise on each finite horizon, so its limits
on overlapping horizons agree and define $u^n$ on $\R_+$.

Set
$M_n(t)=\sup_x\e^{-2\mu|x|}\|u^n(t,x)\|_{L^p(\Omega)}^2$.
By \condH{}(i) and \eqref{eq:ex-mollifier}, there is a constant $C_1>0$
independent of $n$ such that
\[
 \e^{-\mu|y|}\max\!\left\{\|b(v,y,u^n(v,y))\|_{L^p(\Omega)},
                  \|\sigma_n(v,y,u^n(v,y))\|_{L^p(\Omega)}\right\}
 \le C_1(2+M_n(v)^{1/2}).
\]
Then, by the estimates in \eqref{eq:ex-drift-diff} and \eqref{eq:ex-noise-diff},
the mild equation gives
\begin{equation}\label{eq:ex-moment-recursion}
 M_n(t)\le C_3+C_3\int_0^t
   \bigl(1+(t-v)^{-1/2}\bigr)M_n(v)\d v
\end{equation}
for a constant $C_3>0$ independent of $n$. Put
$B_n=\sup_{0\le t\le T}\e^{-\theta t}M_n(t)$, which is finite for
each $n$ by Picard convergence. Multiplying \eqref{eq:ex-moment-recursion} by $\e^{-\theta t}$ and
taking the supremum gives
\[
 B_n\le C_3+C_3(\theta^{-1}+\sqrt{\pi}\,\theta^{-1/2})B_n .
\]
Choose one $\theta>0$ so that the coefficient of $B_n$ is at most
$1/2$. Then \eqref{eq:ex-moment} follows.
\end{proof}

\begin{lemma}[Weak existence]\label{lem:ex-weak}
Assume \condH{}. Then, for every deterministic $u_0\in C_{\rm tem}$,
there is a filtered probability space carrying a space-time white noise
$W$ and a weak solution $u$ with initial function $u_0$ in the sense of
Definition~\ref{def:solution}.
\end{lemma}

\begin{proof}
Let $u^n$ be the approximate solutions from
Lemma~\ref{lem:ex-approx}. The proof has two steps. First,
we use the Kolmogorov--Chentsov theorem for continuous modifications
and its tightness corollary for local tightness
\cite[Theorem~3.23 and Corollary~16.9]{Kallenberg2002}, then use a spatial tail
estimate to obtain tightness in
$C(\R_+;C_{\rm tem})$. Second, we take a joint limit of the approximate
solutions and the Brownian coordinates of the noise. We pass the
martingale and covariation identities to the limit to recover the
Walsh integral in the original equation. This is the approximation
strategy of \cite[proof of Theorem~1.1]{MP}.

\medskip\noindent
\emph{Step 1.}
The estimates below hold for every choice of $\mu>0$;
we may therefore choose the weight parameter separately for each
seminorm of $C_{\rm tem}$.
Fix $p>8$ and $T>0$, and write $Y^n_t=u^n_t-P_tu_0$.
The constant $\widehat C$ below is uniform in $n,j$ and the
space-time variables. It may depend on $p,T,\mu$, $|u_0|_\mu$ and the
coefficient bounds.
By \eqref{eq:ex-moment}, \condH{}(i) and \eqref{eq:ex-mollifier},
\[
 \|b(v,z,u^n_v(z))\|_{L^p(\Omega)}
 +\|\sigma_n(v,z,u^n_v(z))\|_{L^p(\Omega)}
 \le \widehat C\e^{\mu|z|},\qquad v\le T.
\]
For $i=0,1,2$, $q=1,2$, the Gaussian kernel satisfies
\[
 \|\e^{\mu|\cdot|}\partial_x^i p_a\|_q
 \le \widehat C a^{-i/2-(1-1/q)/2},\qquad 0<a\le T.
\]
For $j\in\mathbb Z$, fix $x,y\in[j,j+1]$, put $d=|x-y|$, and take $0\le s<t\le T$
with $\ell=t-s$. The triangle inequality for the exponential
weight and the spatial derivative estimate give
\[
 \|\e^{\mu|\cdot|}(p_a(\cdot+x-y)-p_a)\|_q
 \le \widehat C a^{-(1-1/q)/2}\min\{1,d/\sqrt a\}.
\]
Similarly, $\partial_a p_a=\tfrac12\partial_{xx}p_a$ gives
\[
 \|\e^{\mu|\cdot|}(p_{a+\ell}-p_a)\|_q
 \le \widehat C a^{-(1-1/q)/2}\min\{1,\ell/a\},
 \qquad a+\ell\le T.
\]
The same estimates as in \eqref{eq:ex-drift-diff} and \eqref{eq:ex-noise-diff} give
\begin{align*}
 &\|Y^n_t(x)-Y^n_s(y)\|_{L^p(\Omega)} \\
 &\le C_1\e^{\mu|j|}
 \left[
 d+\ell(2+\log(T/\ell))\right.\\
 &\qquad\left.
 +\left(\int_0^T a^{-1/2}\min\{1,d^2/a\}\d a\right)^{1/2}
 +\left(\int_0^T a^{-1/2}\min\{1,\ell^2/a^2\}\d a
                 +2\sqrt{\ell}\right)^{1/2}
 \right]\\
 &\le C_2\e^{\mu|j|}(d^{1/2}+\ell^{1/4})
\end{align*}
for some constants $C_1,C_2>0$ independent of $n$.
On $Q_j=[0,T]\times[j,j+1]$, the increment estimate reads
\begin{equation}\label{eq:ex-increment}
 \E|Y^n(z)-Y^n(z')|^p
 \le C_3\e^{p\mu|j|}|z-z'|^{p/4},\qquad z,z'\in Q_j,
\end{equation}
for some constant $C_3>0$ independent of $n$ and $j$.
For $0<\alpha<1$, let $C^\alpha(Q_j)$ be the space of functions
$f:Q_j\to\R$ with finite norm
\[
 \|f\|_{C^\alpha(Q_j)}\coloneqq
 \sup_{z\in Q_j}|f(z)|+
 \sup_{\substack{z,z'\in Q_j\\z\ne z'}}
 \frac{|f(z)-f(z')|}{|z-z'|^\alpha},
\]
where $|z-z'|$ is the Euclidean distance.
By the dyadic estimate in the proof of \cite[Theorem~3.23]{Kallenberg2002}
and $Y^n(0,j)=0$,
$(t,x)\mapsto Y^n_t(x)$ has a continuous modification for each $n$, and
\begin{equation}\label{eq:ex-holder-moment}
 \E\|Y^n\|_{C^\alpha(Q_j)}^p\le C_{\alpha}\e^{p\mu|j|},
 \qquad 0<\alpha<\frac14-\frac2p,
\end{equation}
for some constant $C_{\alpha}>0$ independent of $n$ and $j$.
By \cite[Corollary~16.9]{Kallenberg2002},
the laws of $Y^n|_{Q_j}$ are tight in $C(Q_j)$.

For $\lambda>\mu$, let
$H^n_R=\sup_{|j|>R,\,(t,x)\in Q_j}\e^{-\lambda|x|}|Y^n_t(x)|$.
Since $|x|\ge |j|-1$ on $Q_j$, the union bound, Markov's inequality
and \eqref{eq:ex-holder-moment} give
\begin{equation}\label{eq:ex-tail}
\mathbb P(H^n_R>\delta)
 \le \delta^{-p}\sum_{|j|>R}
      \e^{-p\lambda(|j|-1)}\E\|Y^n\|_{C^\alpha(Q_j)}^p
 \le C_\alpha\e^{p\lambda}\delta^{-p}\sum_{|j|>R}\e^{-p(\lambda-\mu)|j|}
 \longrightarrow0
\end{equation}
uniformly in $n$ as $R\to\infty$. For each $n$, $H^n_R$ decreases in $R$, so the above
probability bound implies $H^n_R\downarrow0$ almost surely.
For $t,s\in[0,T]$ and integers $R\ge1$,
\begin{equation}\label{eq:ex-weight-continuity}
 |Y^n_t-Y^n_s|_\lambda
 \le \sup_{|x|\le R+1}|Y^n_t(x)-Y^n_s(x)|+2H^n_R .
\end{equation}
The first term of \eqref{eq:ex-weight-continuity} tends to zero as $s\to t$ by local continuity. Letting
$R\to\infty$ in \eqref{eq:ex-weight-continuity} proves continuity in $|\cdot|_\lambda$.
Taking $T=m$ and $\lambda=1/k$ for $m,k\ge1$, with $\mu<1/k$, gives
$Y^n\in C(\R_+;C_{\rm tem})$ almost surely, and hence the
same for $u^n=Y^n+P_\cdot u_0$.

Local tightness on $K_R=[0,m]\times[-R,R]$, together with \eqref{eq:ex-tail}
and the Arzel\`a--Ascoli theorem \cite[Theorem~7.25]{Rudin1976}, gives tightness of
$\{Y^n\}_{n\ge1}$ on $[0,m]$ in the norm $|\cdot|_{1/k}$.
A diagonal argument over $m,k\in\mathbb N$ then gives
tightness of $(Y^n)$, and hence of $(u^n)$, in the Polish space
$C(\R_+;C_{\rm tem})$.

\medskip\noindent
\emph{Step 2.}
Let $(e_k)_{k\ge1}$ be a real orthonormal basis of $L^2(\R,\d x)$ and
set $B_t^k:=W(1_{[0,t]}e_k)$.
The joint laws of $(u^n,(B^k)_{k\ge1})$ are tight
in $C(\R_+;C_{\rm tem})\times C(\R_+;\R)^{\mathbb N}$, since the law of the
Brownian coordinates is fixed.
By the Skorokhod representation theorem,
after taking a subsequence and changing probability space, we may
assume that
$(u^n,(B^{n,k})_{k\ge1})\to(u,(B^k)_{k\ge1})$
almost surely.
For $\phi\in C_c^\infty(\mathbb R)$, define
\[
 \begin{aligned}
 M_t^\phi
 &:=\langle u_t-u_0,\phi\rangle
 -\int_0^t\left[
 \tfrac12\langle u_v,\phi''\rangle
 +\langle b(v,\cdot,u_v),\phi\rangle
 \right]\d v,\\
 g_v(x)&:=\phi(x)\sigma(v,x,u_v(x)).
 \end{aligned}
\]
Define $M^{n,\phi}$ and $g^n$ by the same formulas with $(u,\sigma)$
replaced by $(u^n,\sigma_n)$.
Write $W^n$ for the white noise represented by $(B^{n,k})_{k\ge1}$ in the filtration of the $n$th pair.
The martingale identities below hold for each represented pair because they depend only on its joint law.
Testing the mild equation for $u^n$ against $\phi$
(see \cite[Theorem~2.1]{Shiga}) gives
\[
 M_t^{n,\phi}
 =\int_0^t\int_{\mathbb R}g_v^n(x)\,W^n(\d v,\d x).
\]
It follows that
\[
 [B^{n,i},B^{n,k}]_t=\mathbf1_{\{i=k\}}t,\qquad
 [M^{n,\phi}]_t=\int_0^t\|g_v^n\|_2^2\,\d v,\qquad
 [M^{n,\phi},B^{n,k}]_t
 =\int_0^t\langle g_v^n,e_k\rangle\,\d v.
\]
The locally uniform convergence $u^n\to u$, the local convergence
$\sigma_n\to\sigma$, and \condH{}(ii)--(iii) imply
$\sup_{t\le T}|M_t^{n,\phi}-M_t^\phi|\to0$
and $\int_0^T\|g_v^n-g_v\|_2^2\,\d v\to0$ as $n\to\infty$ almost surely.
Test increments of the martingales and their compensated products, such as
$M^{n,\phi}_tB^{n,k}_t-\int_0^t\langle g_v^n,e_k\rangle\,\d v$,
against bounded continuous functions of the joint past. At each fixed time,
\eqref{eq:ex-moment} and the Brownian moment bounds give uniform integrability.
Passing to the limit shows that $M^\phi$ and the $B^k$ are continuous
martingales in the joint filtration, with $[B^i,B^k]_t=\mathbf1_{\{i=k\}}t$.
Hence
$[M^\phi]_t=\int_0^t\|g_v\|_2^2\,\d v$ and
$[M^\phi,B^k]_t=\int_0^t\langle g_v,e_k\rangle\,\d v$.
L\'evy's characterization makes the $B^k$ independent Brownian motions
in this filtration. They define a space-time white noise $W$ on the basis
$(e_k)_{k\ge1}$. By the bracket identities and Parseval's identity,
\[
 \E\left|
 M_t^\phi-\sum_{k=1}^l
 \int_0^t\langle g_v,e_k\rangle\,\d B_v^k
 \right|^2=
 \E\int_0^t\left(
 \|g_v\|_2^2
 -\sum_{k=1}^l|\langle g_v,e_k\rangle|^2
 \right)\d v
 \longrightarrow0
\]
as $l\to\infty$. Hence
\[
 M_t^\phi
 =
 \int_0^t\int_{\mathbb R}
 \phi(x)\sigma(v,x,u_v(x))\,W(\d v,\d x).
\]
Thus $u$ is a weak solution to \eqref{eq:spde}.
\end{proof}

\begin{proof}[Proof of Theorem~\ref{thm:main}]
Lemma~\ref{lem:ex-weak} supplies a weak
solution, while Proposition~\ref{prop:stability} gives pathwise uniqueness for
solutions driven by the same white noise in the same filtration.
For each $t$, the noise increments after $t$ are independent of the
solution and the noise up to $t$, so the solution is temporally compatible
with the noise.
Kurtz's Yamada--Watanabe theorem for stochastic models
\cite[Theorem~1.5 and Lemma~2.10]{Kurtz2014} then yields a solution measurable with
respect to the prescribed noise (and hence adapted to its natural
filtration by \cite[Proposition~2.13]{Kurtz2014}), together with uniqueness of the joint law of the solution
and its driving noise.
\end{proof}
\endgroup

\end{document}